\documentclass[11pt]{article}

\usepackage{amsmath,latexsym,epsfig}
\usepackage{amsfonts,amssymb,amscd}

\title{\bf A proof of a conjecture for the  number of ramified coverings of the sphere by the  
torus\footnote{1991 Mathematics Subject Classification: Primary 58D29, 58C35;
Secondary 05C30, 05E05}}

\author{
I.P.Goulden\thanks{Dept. of Combinatorics and Optimization,
University of Waterloo, Waterloo, Ontario, Canada} and
D.M.Jackson\thanks{Dept. of Combinatorics and Optimization,
University of Waterloo, Waterloo, Ontario, Canada} 
}
\date{December 11, 1998}

\begin{document}
 \maketitle

 \newtheorem{theorem}{Theorem}[section]
 \newtheorem{proposition}[theorem]{Proposition}
 \newtheorem{definition}[theorem]{Definition}
 \newtheorem{axiom}[theorem]{Axiom}
 \newtheorem{lemma}[theorem]{Lemma}
 \newtheorem{corollary}[theorem]{Corollary}
 \newtheorem{remark}[theorem]{Remark}
 \newtheorem{example}[theorem]{Example}
 \newtheorem{conjecture}[theorem]{Conjecture}

\def\sP{{\sf{P}}}
\def\sX{{\sf{X}}}
\def\bfp{{\rm\bf p}}
\def\cC{{\cal{C}}}
\def\cH{{\cal{H}}}
\def\cP{{\cal{P}}}
\def\cU{{\cal{U}}}
\def\symgp{{\mathfrak{S}}}
\def\rats{{\mathbb{Q}}}
\def\proof{{\rm\bf Proof:\quad}}
\def\qed{{\hfill{\large $\Box$}}}
\def\sph{\mathbb{S}^2}
\def\surf{\sX}
\def\pdif#1#2{\frac{\partial{#1}}{\partial p_{#2}}}
\def\dpdif#1#2#3{\frac{\partial^2{#1}}{\partial p_{#2}\partial p_{#3}}}
\def\macz#1{\vartheta(#1)}
\def\wij#1#2{w^{(#1)}_{#2}}

\def\atp#1#2{\stackrel{\scriptstyle{#1}}{\scriptstyle{#2}}}



\begin{abstract}
An explicit expression is obtained for the generating series for the
number of ramified coverings of the sphere by the torus, with elementary
branch points and prescribed ramification type over infinity.
This proves a conjecture of Goulden, Jackson and Vainshtein for the explicit number of
such coverings.
\end{abstract}

\section{Introduction}\label{SIn}
Let $\surf$ be a compact connected Riemann surface of genus $g\ge0.$
A {\em ramified covering} of $\sph$ of degree $n$ by $\surf$ is a non-constant
meromorphic function $f\colon\surf\longrightarrow\sph$ such that $\vert f^{-1}(q)\vert=n$
for all but a finite number of points $q\in\sph,$ which are called {\em branch points}.
Two ramified coverings $f_1$ and $f_2$ of $\sph$ by $\surf$ are said to be {\em equivalent}
if there is a homeomorphism $\pi\colon\surf\longrightarrow\surf$ such that $f_1=f_2\pi.$
A ramified covering $f$ is said to be {\em simple} if $\vert f^{-1}(q)\vert=n-1$
for each branch point of $f$, and is {\em almost simple} if $\vert f^{-1}(q)\vert=n-1$
for each branch point but one, that is denoted by $\infty.$
The preimages of $\infty$ are the poles of $f$. If $\alpha_1,\ldots,\alpha_m$ are
the orders of the poles of $f,$ where $\alpha_1\ge\ldots\ge\alpha_m\ge1,$ then
$\alpha=(\alpha_1,\ldots,\alpha_m)$ is a partition of $n$ and is called the
{\em ramification type} of $f.$

Let $\mu_m^{(g)}(\alpha)$ be the number of almost simple ramified coverings
of $\sph$ by $\surf$ with ramification type $\alpha.$ The problem of determining
an (explicit) expression for $\mu_m^{(g)}(\alpha)$ is called the
{\em Hurwitz Enumeration Problem}. The purpose of this paper is to prove the
following result for the torus, giving an explicit expression for $\mu_m^{(1)}(\alpha)$ for
an arbitrary partition $\alpha=(\alpha_1,\ldots,\alpha_m).$
Theorem~\ref{T1} was previously conjectured by Goulden, Jackson and Vainshtein
in~\cite{GJV} where it was proved for all partitions $\alpha$ with $m\le6,$
and for the particular partition $(1^m)$, for any $m\ge1.$
Let $\cC_\alpha$ be the conjugacy class of the symmetric group
$\symgp_n$ on $n$ symbols indexed by the partition $\alpha$ of $n.$

\begin{theorem}\label{T1}
$$
\mu_m^{(1)}(\alpha) = \frac{\vert\cC_\alpha\vert}{24\, n!} (n+m)!
\left( \prod_{i=1}^m\frac{\alpha_i^{\alpha_i}}{(\alpha_i-1)!}\right)
\left(n^m-n^{m-1}-\sum_{i=2}^m (i-2)!e_i n^{m-i}\right)
$$
where $e_i$ is the $i$-th elementary symmetric function in $\alpha_1,\ldots,\alpha_m$
and $e_1=\alpha_1+\cdots+\alpha_m=n.$
\end{theorem}

Previously Hurwitz~\cite{Hrfgv} had shown (see also Goulden and Jackson~\cite{GJtransf} and
Strehl~\cite{vs}) that, for the sphere,
\begin{eqnarray}\label{e0}
\mu_m^{(0)}(\alpha) = \frac{\vert\cC_\alpha\vert}{n!} (n+m-2)! n^{m-3}
\left( \prod_{i=1}^m\frac{\alpha_i^{\alpha_i}}{(\alpha_i-1)!}\right).
\end{eqnarray}
The approach developed by Hurwitz is outlined in the next section.

Very recently Vakil~\cite{V2} has given an independent proof of Theorem~\ref{T1}.
He develops, by techniques in algebraic geometry, and solves a recurrence equation  
that is completely different
in character from the one obtained from the differential equation in this
paper.


\section{Hurwitz's combinatorialization of ramified coverings}
Hurwitz's approach was to represent a ramified covering $f$  of $\sph,$ with ramification type $\alpha,$
by a combinatorial datum $(\sigma_1,\ldots,\sigma_r)$ consisting
of transpositions in $\symgp_n,$ whose product $\pi$ is in $\cC_\alpha,$
such that $\langle \sigma_1,\ldots,\sigma_r\rangle$ acts transitively
on the set $\{1,\ldots,n\}$ of sheet labels and that $r=n+m+2(g-1),$
where $m=l(\alpha),$ the length of $\alpha.$ The latter condition is a
consequence of the Riemann-Hurwitz formula. Under this combinatorialization
he showed that
$$\mu_m^{(g)}(\alpha) =  \frac{\vert\cC_\alpha\vert}{n!} c_g(\alpha),$$
where $c_g(\alpha)$ is the number of  such factorizations of an arbitrary
but fixed $\pi\in\cC_\alpha.$ He studied the action of $\sigma_r$ on
$\sigma_1\cdots\sigma_{r-1}$ to derive a recurrence equation for $c_g(\alpha).$
The difficulty with Hurwitz's approach is that the recurrence equations for $c_g(\alpha)$ are
intractable in all but a small number of special cases.

It appears that his approach can be made more tractable by the introduction
of {\em cut operators} and {\em join operators} that have been developed for combinatorial
purposes by Goulden~\cite{Gdosf}, Goulden and Jackson~\cite{GJtransf} and
Goulden, Jackson and Vainshtein~\cite{GJV}.
These are partial differential operators in indeterminates $p_1, p_2,\ldots$
that take account of the enumerative consequences of the action of $\sigma_r$ on
$\rho=\sigma_1\cdots\sigma_{r-1},$ when summed over all such ordered transitive
factorizations. There are two cases.
The action of $\sigma_r$ on $\rho$ is either to join an $i$-cycle and a $j$-cycle
of $\rho$ to produce an $i+j$-cycle, or to cut an $i+j$-cycle of $\rho$ to produce
an $i$-cycle and a $j$-cycle. In the first case the operators are the join operators
$$
p_{i+j}\frac{\partial^2}{\partial p_i\partial p_j}\quad\mbox{and}\quad
p_{i+j}\left( \frac{\partial}{\partial p_i} \right) \left( \frac{\partial}{\partial p_j} \right),$$
and in the second case the operator is the cut operator
$$p_ip_j \frac{\partial}{\partial p_{i+j}}.$$
A ``cut-and-join'' analysis of the action of $\sigma_r$ on $\rho$ therefore leads
to a nonhomogeneous partial differential equation in a countably infinite
number of variables (indeterminates) for the generating series $\Phi$ for $c_g(\alpha).$
The type of $\Phi$ is determined by the combinatorial properties of the
cut-and-join analysis.

The advantage of this approach to the Hurwitz Enumeration Problem is that
it facilitates the transformation of the differential equation for $\Phi$
by an implicit change of variables. The series that is involved with this
transformation is denoted by $s=s(x,\bfp)$ throughout, where $\bfp=(p_1,p_2,\ldots),$
and appears to be fundamental to the problem. 


\section{The differential equation}\label{Stde}
Let $\bfp_\alpha=p_{\alpha_1}\cdots p_{\alpha_m}$ where $\alpha=(\alpha_1,\ldots,\alpha_m).$
Let
\begin{eqnarray}\label{e1}
\Phi(u,x,z,\bfp) = \sum_{\atp{n,m\ge1}{g\ge0}}
\sum_{\atp{\alpha\vdash n} {l(\alpha)=m}}
\vert\cC_\alpha\vert c_g(\alpha)
\frac{u^{n+m+2(g-1)}}{(n+m+2(g-1))!}
\frac{x^n}{n!} z^g p_\alpha,
\end{eqnarray}
the generating series for $c_g(\alpha),$ where $\alpha\vdash n$ signifies that
$\alpha$ is a partition of $n.$ It was shown in~\cite{GJV} that $f=\Phi(u,1,z,\bfp)$
satisfies the partial differential equation

\begin{eqnarray}\label{e2}
\frac{\partial f}{\partial u} = \frac12 \sum_{i,j\ge1} \left(
ijp_{i+j}z\dpdif fij + ijp_{i+j}
\pdif fi \pdif fj +(i+j)p_ip_j \pdif f{i+j} \right).
\end{eqnarray}
By replacing $p_i$ by $x^ip_i$ for $i\ge1$ it is readily seen that
$f=\Phi(u,x,z,\bfp)$ satisfies~(\ref{e2}). But 
$\Phi(u,x,z,\bfp)\in\rats[u,z,\bfp]\,[[x]],$
and it is also readily seen that~(\ref{e2}) has a unique solution in this ring.

Let $F_i(x,\bfp)=[z^i]\Phi(1,x,z,\bfp)$ for $i=0,1,$ 
where $[z^i]f$ denotes the coefficient of $z^i$ in the formal power series $f.$
Then $F_0$ is the
generating series for the numbers $c_0(\alpha)$, which have been determined by Hurwitz, so $F_0$ is
known. $F_1$ is the generating series for $c_1(\alpha).$
The next result gives the linear first order partial differential equation for $F_1$
that is induced by restricting~(\ref{e2}) above to terms of degree at most one in $z$.

\begin{lemma}\label{L1}
The series $f=F_1$ satisfies the partial differential equation
\begin{eqnarray}\label{e3}
T_0 f- T_1=0
\end{eqnarray}
where
\begin{eqnarray*}
T_0 &=& x\frac{\partial}{\partial x} + \sum_{i\ge1} p_i\pdif {}i
- \sum_{i,j\ge1} ij p_{i+j} \pdif{F_0}i \pdif{}j
-\frac12 \sum_{i,j\ge1} (i+j) p_ip_j \pdif{}{i+j}, \\
T_1 &=& \frac12 \sum_{i,j\ge1} ij p_{i+j} \dpdif{F_0}ij.
\end{eqnarray*}
\end{lemma}
\proof Clearly, from~(\ref{e2}),
$$u\frac{\partial}{\partial u} [z] \Phi = 
[z] \left( 
x\frac{\partial}{\partial x} + \sum_{i\ge1} p_i\pdif {}i
\right) 
\Phi.$$
The result follows by applying $[z]$ to~(\ref{e2}). \qed

We now turn our attention to solving this partial differential equation.
Let
\begin{eqnarray}\label{e4}
G_1(x,\bfp) = \frac{1}{24}
\sum_{n,m\ge1}
\sum_{\atp{\alpha\vdash n} {l(\alpha)=m}}
\vert\cC_\alpha\vert 
\left( \prod_{i=1}^m\frac{\alpha_i^{\alpha_i}}{(\alpha_i-1)!}\right)
\left(n^m-n^{m-1}-\sum_{i=2}^m (i-2)!e_i n^{m-i}\right)
\frac{x^n}{n!}p_\alpha.
\end{eqnarray}
Since~(\ref{e2}) has a unique solution in $\rats[u,z,\bfp]\,[[x]],$ then~(\ref{e3}) has
a unique solution in  $\rats[\bfp]\,[[x]].$ To establish Theorem~\ref{T1}
it therefore suffices to show that $f=G_1$ satisfies~(\ref{e3})
(note that $G_1$ has a constant term of $0,$ so the initial condition is satisfied).


\section{The generating series $G_1$}\label{StgsG}
To obtain a convenient form for $G_1$ the following lemma is required that 
expresses the elementary symmetric function $e_k(\lambda)$ as  the coefficient
in a formal power series. For a partition 
$\alpha=(\alpha_1,\ldots,\alpha_r)$ let $m_i$ denote the number of occurrences
of $i$ in $\alpha,$ and  we may therefore write $\alpha=(1^{m_1},\ldots,r^{m_r}).$ 
Let 
$\macz\alpha=\prod_{i=1}^r i^{m_i} m_i!.$ Let $\cP$ denote the set of all
partitions with the null partition adjoined.

\begin{lemma}\label{Lvart}
For any nonnegative integer $k$ and partition $\lambda,$
$$e_k(\lambda) = 
\frac{\macz\lambda}{k!}[p_\lambda](p_1+p_2+\cdots)^k
\sum_{\alpha\in\cP} \frac{p_\alpha}{\macz\alpha}.$$
\end{lemma}
\proof First
$$\sum_{\alpha\in\cP} \frac{p_\alpha}{\macz\alpha} = 
\sum_{m_1,m_2,\ldots\ge0}
\frac{p_1^{m_1}}{1^{m_1}m_1!}\frac{p_2^{m_2}}{2^{m_2}m_2!}\cdots$$
and
$$(p_1+p_2+\cdots)^k = \sum_{\atp{i_1,i_2,\ldots\ge0} {i_1+i_2+\cdots=k}}
k!\frac{p_1^{i_1}p_2^{i_2}\cdots}{i_1!i_2!\cdots}.$$
If  $\macz\lambda=1^{j_1}j_1!2^{j_2}j_2!\cdots$ then
\begin{eqnarray*}
\frac{\macz\lambda}{k!}[p_\lambda](p_1+p_2+\cdots)^k
\sum_{\alpha\in\cP} \frac{p_\alpha}{\macz\alpha} &=&
\sum_{\atp{m_1,i_1,m_2,i_2,\ldots\ge0} {i_1+i_2+\cdots=k}}
\binom{j_1}{i_1}1^{i_1} \binom{j_2}{i_2}2^{i_2} \cdots
\end{eqnarray*}
where the sum is further  constrained by $i_1+m_1=j_1, i_2+m_2=j_2,\ldots.$
Then
\begin{eqnarray*}
\frac{\macz\lambda}{k!}[p_\lambda](p_1+p_2+\cdots)^k
\sum_{\alpha\in\cP} \frac{p_\alpha}{\macz\alpha} &=&
[t^k](1+t)^{j_1}(1+2t)^{j_2}(1+3t)^{j_3}\cdots \\
&=& e_k( (1^{j_1},2^{j_2},\ldots) )
\end{eqnarray*}
and the result follows. \qed

Let
\begin{eqnarray}\label{e5}
\psi_i(x,\bfp)=\sum_{r\ge1}r^{i-1}a_rp_r x^r,
\end{eqnarray}
where $i$ is an integer and $a_r=r^r/(r-1)!,$ for $r\ge1.$ Let
$s\equiv s(x,\bfp)$ be the unique solution of the functional equation
\begin{eqnarray}\label{e6}
s = x e^{\psi_0(s,\bfp)}
\end{eqnarray}
in the ring $\rats[\bfp]\,[[x]].$ An explicit series expansion of $s$ can be obtained  by
Lagrange's Implicit Function Theorem (see~\cite{GJbook}, for example).
Let $\psi_i$ denote $\psi_i(s,\bfp).$ The next result
gives a closed form expression for $G_1$ in terms of $s,$ and indicates
the fundamental importance of the series $s$ to the solution of the partial
differential equation given in~(\ref{e2}).

\begin{theorem}\label{T2}
$$G_1(x,\bfp)=\frac{1}{24}\log(1-\psi_1)^{-1}-\frac{1}{24}\psi_0.$$
\end{theorem}
\proof For a partition $\alpha=(\alpha_1,\ldots,\alpha_m),$ let
$a_\alpha=a_{\alpha_1}\cdots a_{\alpha_m}$ and
$$g(x,\bfp)=\sum_{n\ge1}\sum_{\atp{\alpha\vdash n} {l(\alpha)=m}}
\frac{n^{m-1}}{\macz\alpha }a_\alpha p_\alpha x^n,$$
a constituent of the series $G_1$ given in~(\ref{e4}), since $\macz\alpha=n!/\vert\cC_\alpha\vert.$
It is easily shown that
$$x\frac{\partial g}{\partial x} = \sum_{n\ge1} x^n[t^n] e^{n\psi_0(t,\bfp)}$$
so, by Lagrange's Implicit Function Theorem,
$$x\frac{\partial g}{\partial x} =  \frac{\psi_1}{1-\psi_1}.$$
But, from~(\ref{e6}),
\begin{eqnarray}\label{e7}
x\frac{\partial s}{\partial x} =  \frac{s}{1-\psi_1}
\end{eqnarray}
and  from~(\ref{e5}),
\begin{eqnarray*}
\frac{\partial \psi_i}{\partial s} =  \frac{1}{s}\psi_{i+1}.
\end{eqnarray*}
Then $x\partial(g-\psi_0)/\partial x = 0$ and, since $g(0,\bfp)=0,$ it follows
that $g(x,\bfp)=\psi_0.$

Next we consider the terms of $G_1$ in~(\ref{e4}) that are not included in
$g(x,\bfp).$ First, note that
$$\sum_{\theta\in\cP}\frac{p_\theta}{\macz\theta} = 
\exp\sum_{i\ge1}\frac{p_i}i,$$
so, replacing $p_i$ by $ntp_i a_i$ for $i\ge1$ in Lemma~\ref{Lvart}, we have
$$
\frac{\macz\alpha}{k!}[p_\alpha t^n]
\left( e^{\psi_0(t,\bfp)}\right)^n \psi_1^k(t,\bfp)
=n^{m-k}a_\alpha e_k(\alpha),
$$
where $m=l(\alpha).$
Then
\begin{eqnarray*}
a_\alpha n^m-\sum_{k\ge2}(k-2)! n^{m-k}a_\alpha e_k(\alpha) &=&
\macz\alpha [p_\alpha t^n]
\left(1-\sum_{k\ge2}\frac1{k(k-1)}\psi_1^k(t,\bfp)\right) \left( e^{\psi_0(t,\bfp)}\right)^n \\
&=&  \macz\alpha [p_\alpha x^n]
\frac1{1-\psi_1} \left(1-\sum_{k\ge2}\frac1{k(k-1)}\psi_1^k\right)
\end{eqnarray*}
by Lagrange's Implicit Function Theorem, for $n\ge1.$
Then
$$a_\alpha n^m-\sum_{k\ge2}(k-2)! n^{m-k}a_\alpha e_k(\alpha) = \macz\alpha [p_\alpha x^n]
\left(1+\log(1-\psi_1)^{-1}\right).$$
The result follows by combining the two expressions that have been obtained
and by using the fact that $G_1(0,\bfp)=0.$ \qed


\section{The proof of Theorem~\ref{T1}}\label{Stfs}
The remaining portion of the paper is concerned with the proof of Theorem~\ref{T1}.


\noindent \proof
Our strategy is to show that the expression for 
$G_1$ given in Theorem~\ref{T2} satisfies the partial differential
equation given in~(\ref{e3}).
We begin by considering  the derivatives that are required in the determination of $T_0G_1-T_1.$
From the functional equation~(\ref{e6})
\begin{eqnarray}\label{e8}
\pdif sk =\frac1k \frac{a_k s^{k+1}}{1-\psi_1}.
\end{eqnarray}
Then, for $k\ge1,$
\begin{eqnarray}\label{e9}
\pdif {\psi_j}k =  k^{j-1}a_ks^k + \frac{a_k}{k} \frac{\psi_{j+1}s^k}{1-\psi_1}.
\end{eqnarray}
The only derivatives of $F_0$ that are needed are
\begin{eqnarray}\label{e10}
\pdif {F_0}k = \frac{a_k}{k^3}s^k -\frac{a_k}{k^2}\sum_{r\ge1}  a_r p_r\frac{s^{k+r}}{k+r},
\end{eqnarray}
from Proposition~3.1~\cite{GJtransf} and, from~(\ref{e8}) and~(\ref{e10}),
\begin{eqnarray}\label{e11}
\dpdif {F_0}ij = \frac{a_ia_j}{ij} \frac{s^{i+j} }{i+j},
\end{eqnarray}
for $i,j\ge1.$
For completeness we note that, from Proposition~3.1~\cite{GJtransf},
$$\left(x\frac{\partial}{\partial x}\right)^2F_0 = \psi_0.$$
The derivatives of $G_1$ that are needed are, from~(\ref{e8}),
\begin{eqnarray}\label{e12}
x\frac{\partial G_1}{\partial x} 
= \frac1{24}\left(\frac{\psi_2}{(1-\psi_1)^2}-\frac{\psi_1}{1-\psi_1}\right)
\end{eqnarray}
and, from~(\ref{e9}), for $k\ge1,$
\begin{eqnarray}\label{e13}
\pdif{G_1}k = \frac{1}{24}a_k\frac{s^k}{1-\psi_1}
+\frac1{24} \frac{a_k}{k}s^k \left( \frac{\psi_2}{(1-\psi_1)^2}-\frac{1}{1-\psi_1} \right).
\end{eqnarray}
Then from Lemma~\ref{L1} and expressions~(\ref{e10}),~(\ref{e11}),~(\ref{e12}) and~(\ref{e13}) it follows that 

\begin{eqnarray}
24(1-\psi_1)^2(T_0G_1-T_1) 
&=&\psi_2(1+\psi_0)-\psi_0(1-\psi_1)-12(1-\psi_1)^2A - (1-\psi_1)B \nonumber \\
&\mbox{}& +(1-\psi_1)C -(\psi_1+\psi_2-1)D+(\psi_1+\psi_2-1)E, \label{e14}
\end{eqnarray}
where
\begin{eqnarray*}
A &=& \sum_{i,j\ge1} \frac{a_ia_j}{i+j}p_{i+j}s^{i+j}, \\
B &=& \sum_{i,j\ge1} \frac{ia_ia_j}{j^2}p_{i+j}s^{i+j}, \\
C &=& \sum_{i,j,m\ge1} \frac{ia_ia_ja_m}{j(j+m)}p_{i+j}p_ms^{i+j+m}
       -\frac12\sum_{i,j\ge1} (i+j)a_{i+j}p_ip_js^{i+j}, \\
D &=& \sum_{i,j\ge1} \frac{a_ia_j}{j^2}p_{i+j}s^{i+j}, \\
E &=& \sum_{i,j,m\ge1} \frac{a_ia_ja_m}{j(j+m)}p_{i+j}p_ms^{i+j+m}
       -\frac12\sum_{i,j\ge1} a_{i+j}p_ip_js^{i+j}.
\end{eqnarray*}

When the expression~(\ref{e14})  is transformed by replacing
$p_is^i$ by $y_i$, for $i\ge1,$ it is  immediately seen to be a polynomial
in $y_1,y_2,\ldots$ of degree $3$ with rational coefficients. If $U_i$ denotes
the degree $i$ part of the expression, for $i=1,\ldots,3,$ the transformed expression
can be written in the form

$$24(1-\psi_1)^2 (T_0G_1-T_1) = U_1+U_2+U_3$$
where
\begin{eqnarray*}
U_1 &=& \psi_2-\psi_0-12A-B+D, \\
U_2 &=& \psi_0(\psi_1+\psi_2)+24\psi_1A+\psi_1B+C-(\psi_1+\psi_2)D-E,\\
U_3 &=& -12\psi_1^2A-\psi_1C+(\psi_1+\psi_2)E.
\end{eqnarray*}
Then $U_i\in\cH_i[y_1,y_2,\ldots],$ the set of homogeneous polynomials of degree $i$
in $y_1, y_2,\ldots.$ Let 
$$\varpi_{1,\ldots,i}\colon\cH_i[y_1,y_2,\ldots]\longmapsto\rats[x_1,x_2,\ldots]$$
be  the symmetrization operation defined by
$$\varpi_{1,\ldots,i} (y_{\alpha_1}\cdots y_{\alpha_i})=\sum_{\pi\in\symgp_i}
x^{\alpha_1}_{\pi(1)}\cdots x^{\alpha_i}_{\pi(i)},$$
extended linearly to $\cH_i[y_1,y_2,\ldots].$ Then $\varpi_{1,\ldots,i}f=0$ implies that $f=0$
for $f\in\cH_i[y_1,y_2,\ldots].$ 

We therefore prove that $U_i=0$ by proving that
$\varpi_{1,\ldots,i}U_i=0,$ for $i=1,2,3.$ To determine the action of the symmetrization operator
on $A,\ldots,E$ it is convenient to introduce the series
$w=w(x)$ as the unique solution of the functional equation
\begin{eqnarray}\label{b1}
w=xe^w
\end{eqnarray}
in the ring $\rats[[x]].$ By Lagrange's Implicit Function Theorem we have
$$w=\sum_{n\ge1}\frac{n^{n-1}}{n!}x^n.$$
Now let $w_i=w(x_i)$ and $w_i^{(j)}=(x_i\partial/\partial x_i)^j w_i.$
Then, from~(\ref{b1}),
\begin{eqnarray}\label{b2}
\wij1{i} = \frac {w_i}{1-w_i},\quad
\wij2{i} = \frac {w_i}{(1-w_i)^3}, \quad
\wij3{i} = \frac {w_i+2w_i^2}{(1-w_i)^5}.
\end{eqnarray}
The action of the symmetrizing operator on $A,\ldots,E$ and their products with $\psi_i$
can be determined in terms of these as follows.

It is readily seen that
\begin{eqnarray*}
\varpi_1(\psi_m) &=& \wij{m+1}1, m\ge-1.
\end{eqnarray*}

For $\varpi_{1}(A),$ using~(\ref{b2}), we have
\begin{eqnarray*}
\varpi_{1}(A) = \sum_{k\ge1} \frac{x_1^k}{k}[x_1^k] \left(w_1^{(2)}\right)^2 =
\int_0^{x_1}\left(w_1^{(2)}\right)^2\frac{dx_1}{x_1}
=\int_0^{w_1} \frac{w_1}{(1-w_1)^5}dw_1
\end{eqnarray*}
so,  by rearrangement
\begin{eqnarray*}
\varpi_1(A) &=& \frac{1}{12}\left((1-w_1)\wij31+w_1\wij21-\wij11\right). 
\end{eqnarray*}

Trivially, 
$$\varpi_1(B) = \wij31 w_1.$$

Next, $\varpi_{1,2}(C)$ is the symmetrization of
$$\sum_{i,j,m\ge1} \frac{ia_ia_ja_m}{j(j+m)}x_1^{i+j}x_2^m
       -\frac12\sum_{i,j\ge1} (i+j)a_{i+j}x_1^ix_2^j$$
with respect to $x_1$ and $x_2.$ Now
\begin{eqnarray*}
\sum_{i,j,m\ge1} \frac{ia_ia_ja_m}{j(j+m)}x_1^{i+j}x_2^m
&=& w_1^{(3)} \sum_{m\ge1} a_mx_2^m \sum_{j\ge1}\frac{a_j}{j}\frac{x_1^j}{j+m} \\
&=& w_1^{(3)} \sum_{m\ge1} a_mx_2^m \frac{1}{x_1^m}\int_0^{x_1} w_1^{(1)}x_1^{m-1}dx_1.
\end{eqnarray*}
But, from~(\ref{b2}), and integrating by parts, we obtain
$$\int_0^{x_1} w_1^{(1)}x_1^{m-1}dx_1
=\int_0^{w_1} w_1^m e^{-mw_1}dw_1
=\frac{1}{a_m}\left( 1-x_1^m \sum_{i=1}^m \frac{m^{m-i}}{(m-i)!}\frac{1}{w_1^i} \right) -\frac{x_1^m}{m}.
$$
Thus
\begin{eqnarray*}
\sum_{i,j,m\ge1} \frac{ia_ia_ja_m}{j(j+m)}x_1^{i+j}x_2^m
&=& w_1^{(3)} \left( \frac{x_2}{x_1-x_2} - \sum_{m\ge1}  x_2^m \sum_{i=1}^m \frac{m^{m-i}}{(m-i)!}\frac{1}{w_1^i} 
-w_2^{(1)} \right) \\
&=& w_1^{(3)} \left( \frac{x_2}{x_1-x_2} - \sum_{m\ge1}x_2^m[t^m] e^{mt} 
\left(\left(1-\frac{t}{w_1}\right)^{-1} -1\right)
-w_2^{(1)} \right) \\
&=&  w_1^{(3)} \left( \frac{x_2}{x_1-x_2} - 
\frac{w_2}{w_1-w_2}\,\frac{1}{1-w_2}
-w_2^{(1)} \right)
\end{eqnarray*}
by the Lagrange Implicit Function Theorem.
Moreover, it is easily seen that
$$
\sum_{i,j\ge1} (i+j)a_{i+j}x_1^ix_2^j = \frac{x_2w_1^{(3)}-x_1w_2^{(3)}}{x_1-x_2}.
$$
Thus, by symmetrizing  the indicated linear combination of these sums, we have
\begin{eqnarray*}
\varpi_{1,2}(C) &=& -\wij31\wij12 - \wij11\wij32
-\frac{\wij31\wij12-\wij11\wij32}{w_1-w_2}.
\end{eqnarray*}

Trivially,
$$\varpi_1(D) = \wij21 w_1.$$

Finally, $\varpi_{1,2}(E)$ is obtained in a  fashion similar to $\varpi_{1,2}(C)$. The expression is
\begin{eqnarray*}
\varpi_{1,2}(E) &=& -\wij21\wij12 - \wij11\wij22
-\frac{\wij21\wij12-\wij11\wij22}{w_1-w_2}.
\end{eqnarray*}

These results may be combined to give expressions for the symmetrizations of $U_1, U_2,U_3$ as follows.

For the  term of degree one,
\begin{eqnarray*}
\varpi_1(U_1) &=&
\wij31-\wij11  -\left( (1-w_1)\wij31+w_1\wij21-\wij11\right)
-\wij31w_1+\wij21w_1.
\end{eqnarray*}
For the term of degree two, after rearrangement,
\begin{eqnarray*}
\varpi_{1,2}(U_2) &=& \left( \wij21\wij32+\wij31\wij22\right)(2-w_1-w_2)
+\wij21\wij22(w_1+w_2)\\
&\mbox{}&  -\frac{\wij31\wij12-\wij11\wij32}{w_1-w_2} + \frac{\wij21\wij12-\wij11\wij22}{w_1-w_2}.
\end{eqnarray*}
When multiplied by $w_1-w_2$ and a suitable power of  $(1-w_1)^{-1}$ and  $(1-w_2)^{-1}$
this becomes a polynomial in $w_1$ and $w_2$ that is identically zero.

For the term of degree three, after rearrangement,
\begin{eqnarray*}
\varpi_{1,2,3}(U_3) 
&=&\frac{1}{w_2-w_3}\left( \wij21\left(\wij32\wij13-\wij12\wij33\right)
-\left(\wij21+\wij31\right)\left(\wij22\wij13-\wij12\wij23\right)\right) \\
&\mbox{}&
+\frac{1}{w_1-w_3}\left( \wij22\left(\wij31\wij13-\wij11\wij33\right)
-\left(\wij22+\wij32\right)\left(\wij21\wij13-\wij11\wij23\right)\right) \\
&\mbox{}&
+\frac{1}{w_1-w_2}\left( \wij23\left(\wij31\wij12-\wij11\wij32\right)
-\left(\wij23+\wij33\right)\left(\wij21\wij12-\wij11\wij22\right)\right) \\
&\mbox{}& -2\wij21\wij22\wij33(1-w_3) - 2\wij21\wij23\wij32(1-w_2) - 2\wij22\wij23\wij31(1-w_1) \\
&\mbox{}& -2 \wij21\wij22\wij23(w_1+w_2+w_3).
\end{eqnarray*}
When multiplied by $(w_1-w_2)(w_2-w_3)(w_1-w_3)$ and a suitable power of  
$(1-w_1)^{-1},$   $(1-w_2)^{-1}$ and  $(1-w_3)^{-1}$
this becomes a polynomial in $w_1$, $w_2$ and $w_3$ that is identically zero.

It is  quickly seen that $\varpi_{1}(U_1)$ is zero. For both $\varpi_{1,2}(\cU_2)$
and $\varpi_{1,2,3}(\cU_3),$ however, the polynomial
expressions were sufficiently large that it was convenient to use {\sf Maple} to
carry out the routine simplification of this stage.

Thus the symmetrization of $24(1-\psi_1)^2(T_0G_1-T_1)=0$ so
$T_0G_1-T_1=0$. It follows from Lemma~\ref{L1} that $F_1=G_1$ and this completes the proof
of Theorem~\ref{T1}.
\qed

\section*{Acknowledgements}
This work was supported by grants individually
to IPG and DMJ from the Natural
Sciences and Engineering Research Council of Canada.




\begin{thebibliography}{999}
 
\bibitem{Gdosf}{\sc I.P.Goulden},
{\em  A differential operator  for symmetric functions and the combinatorics
of multiplying  transpositions},
Trans. Amer. Math. Soc.,
{\bf344} (1994), 421--440.

\bibitem{GJtransf}{\sc I.P.Goulden and D.M.Jackson},
{\em Transitive factorizations into transpositions and holomorphic mappings on the sphere},
Proc. Amer. Math. Soc.,
{\bf125} (1997), 51--60.

\bibitem{GJbook}{\sc I.P.Goulden and D.M.Jackson}
``Combinatorial Enumeration,''
Wiley, New York,
1983.

\bibitem{GJV}{\sc  I.P.Goulden, D.M.Jackson and A.Vainshtein},
{\em The number of ramified coverings of the sphere by the torus
and surfaces of higher genera},
(preprint) September 1998.

\bibitem{Hrfgv}{\sc A.Hurwitz},
{\em Ueber Riemann'sche Fl\"{a}chen mit gegebenen Verzweigungspunkten},
Matematische Annalen,
{\bf39} (1891), 1--60.

\bibitem{vs}{\sc V.Strehl},
{\em  Minimal  transitive products  of transpositions - the reconstruction of a proof by A.Hurwitz},
S\'{e}m. Lothar. Combinat.  {\bf37} (1996), Art.~S37c, pp.~12.

\bibitem{V2}{\sc R.Vakil},
{\em Recursions, formulas, and graph-theoretic interpretations
of ramified coverings of the sphere by surfaces of genus 0 and 1},
CO/9812105.

\end{thebibliography}
\end{document}